\documentclass{article}
\usepackage{amsthm, amssymb, amsfonts}
\author{{\normalsize Samsonov M.E.}\footnote{E-mail address: 
samsonov@pink.phys.spbu.ru}
\\[1ex]{\normalsize Theoretical department}\\
{\normalsize St. Petersburg University}\\
{\normalsize Institute of Physics}}
\title{Quantization of semi-classical twists and noncommutative geometry}
\date{}
\begin{document}
\theoremstyle{plain}\newtheorem{aaa}{Proposition}\newtheorem{ccc}{Theorem}
\theoremstyle{remark}\newtheorem{bbb}{\rm\bf Example}
\theoremstyle{definition}\newtheorem{vvv}{Definition}
\maketitle
\begin{abstract}
A problem of defining the quantum analogues for semi-classical 
twists in $U(\mathfrak{g})[[t]]$ is considered. First, we study 
specialization at $q=1$ of singular coboundary twists defined 
in $U_{q}(\mathfrak{g})[[t]]$ for $\mathfrak{g}$ 
being a nonexceptional Lie algebra, then we consider specialization 
of noncoboundary twists when $\mathfrak{g}=\mathfrak{sl}_{3}$ and 
obtain $q-$deformation of the semi-classical twist introduced by Connes 
and Moscovici in noncommutative geometry.        
\end{abstract}
\hspace{0.9 cm}{\it Keywords:} Noncommutative geometry, Hopf algebras
\section{Introduction}
Hopf algebras play an increasingly important role in noncomutative 
geometry and Quantum Field Theory. One of the sources for producing 
new types of Hopf algebras is twisting, a deformation of the 
coalgebraic structure of a given Hopf algebra  
$(H,\mu,\eta,\Delta,\epsilon,S)$ preserving the algebraic structure 
$(H,\mu,\eta)$. Such deformations 
are generated by the twisting elements (twists) 
${\cal F}\in (H\otimes H)$ satisfying the conditions 
\begin{equation}
{\begin{array}{ccc}
{\cal F}^{12}(\Delta\otimes{\rm id})({\cal F})&=&{\cal F}^{23}({\rm 
id}\otimes\Delta)({\cal F})\\[2ex]
(\varepsilon\otimes{\rm id})({\cal F})&=&({\rm id}\otimes\varepsilon)
({\cal F})=1
\end{array}}
\label{drinfeld}
\end{equation}
that guarantee 
$H^{\cal F}\equiv(H,\mu,\eta,{\rm Ad}{\cal F}\circ\Delta,\epsilon,S)$ 
is a new Hopf algebra. In fact, when $H$ is not finite dimensional,
${\cal F}$ is usually defined in some completion of the tensor product and 
$H$ is understood to be a topological Hopf algebra.

In this article we consider two types of twists: 
the semi-classical ones if $H=~U(\mathfrak{g})[[t]]$ and the quantum ones 
if $H=U_{q}(\mathfrak{g})[[t]]$. Some of the semi-classical deformations 
such as those defined by the Jordanian twists \cite{GZ,KLM} appear as the limiting cases 
of the quantum ones (in the sense that specialization at $q=1$ is extended 
to work for topological Hopf algebras). 
It was a motivation for us to study the quantum twists as many 
computational problems involved into a direct check of (\ref{drinfeld}) 
drastically  resolve when one works with $U_{q}(\mathfrak{g})[[t]]$ instead 
of $U(\mathfrak{g})[[t]]$ and thus the quantum twists is a source for many 
universal deformation formulas in the sense of \cite{GZ}.

The work is organized as follows. After preliminary section intended to 
fix notations, we show that if $\mathfrak{g}$ is a nonexceptional simple Lie 
algebra, then a quantum analogue of the Jordanian twist can be taken to be a 
coboundary twist in 
$U_{q}(\mathfrak{g})[[t]]$:
$$
{\cal J}(e_{\lambda}):=(W\otimes W)\Delta(W^{-1}), \mbox{ where }
W=\exp_{q_{\lambda}}(\frac{t}{1-q_{\lambda}}\mathop{}e_{\lambda}); 
\hphantom{aaa}q_{\lambda}:=q^{(\lambda,\lambda)}  
$$
with $e_{\lambda}$ being a quantum highest root generator in some 
quantum Cartan-Weyl basis. We prove that 
${\cal J}(e_{\lambda})$ is nonsingular and specializes to a nontrivial 
twisting of $U(\mathfrak{g})[[t]]$. As an application of the Jordanian twists 
\cite{GZ,KLM} to noncommutative geometry \cite{CM}, we prove that there is a 
homomorphism of the Connes-Moscovici Hopf algebra 
$\iota:{\cal H}_{1}\rightarrow U^{\cal F}(\mathfrak{sl}_{3})[[t]]$, 
with ${\cal F}$ being a Jordanian twist, where ${\cal H}_{1}$ has the 
following structure   
$$
\begin{array}{lclclclcl} 
[Y,X]=X,&& [Y,\delta_{n}]=n\delta_{n},&&
[X,\delta_{n}]=\delta_{n+1},&&[\delta_{k},\delta_{l}]=0,&&k,l\ge 1
\end{array}
$$ 
\begin{equation}
\label{structure}
\begin{array}{c}
{ 
\begin{array}{lcl}
\Delta(Y)=Y\otimes 1+1\otimes Y,&& \Delta(\delta_{1})=
\delta_{1}\otimes 1+1\otimes\delta_{1}\\[2ex]
\end{array}
}\\[2ex]
\Delta(X)=X\otimes 1+1\otimes X+\delta_{1}\otimes Y.
\end{array}
\label{coproducts}
\end{equation}

Through factoring ${\cal H}_{1}^{\prime}:={\cal H}_{1}/<\delta_{2}-\frac 12\mathop{}\delta_{1}^{2}>$,
one obtains in fact an embedding  
$$
\iota: {\cal H}_{1}^{\prime}\hookrightarrow 
U^{\cal F}(\mathfrak{sl}_{3})[[t]]
$$  
and the twist found in \cite{CM}:
\begin{equation}
F=\sum_{n\ge 0}t^{n}\sum_{k=0}^{n}
\frac{S(X)^{k}}{k!}(2Y+k)_{n-k}\otimes\frac{X^{n-k}}{(n-k)!}(2Y+n-k)_{k}
\label{connestwist}
\end{equation}
where
$(\alpha)_{k}=\alpha(\alpha+1)\cdots (\alpha+k-1)$ and 
$S(X)=-X+\delta_{1}Y$, can be obtained as a pullback $F=\iota_{*}\Phi$ of 
a semi-classical twist $\Phi$ in $U^{\cal F}(\mathfrak{sl}_{3})[[t]]$. 
In section \ref{connes} we show that $\iota$ can be "quantized", thus leading 
to a quantum analogue of ${\cal H}_{1}^{\prime}$: 
$$
\begin{array}{lclcl}
k xk^{-1}=q^{2}\mathop{}x,&& k z k^{-1}=q^{2}\mathop{}z,
&&
q^{2}xz-zx=-t z^{2}\\[2ex]
\end{array}
$$ 
$$
\begin{array}{lcl}
\Delta(k)=k\otimes k, &&
\Delta(z)=z\otimes k+1\otimes z
\displaystyle
\end{array}
$$
$$
\Delta(x)=x\otimes k^{-1}+1\otimes x+t\mathop{}z\otimes 
\frac{(k-k^{-1})}{1-q^{2}}.
$$
\subsection*{Acknowledgment}
I would like to thank V. Lyakhovsky, A. Stolin and V. Tolstoy for
valuable comments on the subject
\section{Preliminarities}
Let $\mathfrak{g}$ be a simple Lie algebra with the set of simple roots 
$\pi=~\{\alpha_{1},\ldots,\alpha_{N}\}$ and the Cartan matrix 
$ (A)_{ij}=a_{ij}=2(\alpha_{i},\alpha_{j})/(\alpha_{i},\alpha_{i})$.
By definition, a Hopf algebra $U_{q}(\mathfrak{g})$ is generated by
$\{e_{i},f_{i}, k_{i}^{\pm 1}\}_{1\le i\le N}$ over $\mathbb{C}(q)$  
which are subject to the following relations
\begin{equation}
k_{i}e_{j}k^{-1}_{i}=q^{(\alpha_{i},\alpha_{j})}e_{j},
\hphantom{aaa}
k_{i}f_{j}k_{i}^{-1}=q^{-(\alpha_{i},\alpha_{j})}f_{j}
\label{rel1}
\end{equation}
\begin{equation}
e_{i}f_{j}-f_{j}e_{i}=\delta_{ij}\displaystyle\frac 
{k_{i}-k_{i}^{-1}}{q-q^{-1}}
\label{rel3}
\end{equation}
\begin{equation}
\sum_{n=0}^{1-a_{ij}}(-1)^{n}\left [1-a_{ij}\atop 
n\right]_{q_{i}}(e_{i})^{n}
e_{j}(e_{i})^{1-a_{ij}-n}=0\hphantom{a}\mbox{for $i\ne 
j$}
\label{rel3}
\end{equation}
\begin{equation}
\sum_{n=0}^{1-a_{ij}}(-1)^{n}\left [1-a_{ij}\atop 
n\right]_{q_{i}}(f_{i})^{n}
f_{j}(f_{i})^{1-a_{ij}-n}=0\hphantom{a}\mbox{for $i\ne 
j$}
\label{rel4}
\end{equation}
where $q_{i}=q^{\frac{(\alpha_{i},\alpha_{i})}{2}}$ and  
$$
\left[m\atop n\right]_{q}\equiv\frac {[m]_{q}!}{[n]_{q}![m-n]_{q}!}
$$
$$
[k]_{q}!\equiv [1]_{q}[2]_{q}\ldots [k]_{q},\hphantom{aa}
[l]_{q}\equiv (q^{l}-q^{-l})/(q-q^{-1})
$$
The Hopf algebra structure is defined uniquely by fixing 
the values of the coproduct on the Chevalley generators
\begin{equation}
\Delta(k_{i})=k_{i}\otimes k_{i}
\label{rel6}
\end{equation}
\begin{equation}
\Delta(e_{i})=k_{i}^{-1}\otimes e_{i}+e_{i}\otimes 1,\hphantom{aa} 
\Delta(f_{i})=f_{i}\otimes k_{i}+1\otimes f_{i}
\label{rel7}
\end{equation}
\begin{eqnarray}
S(k_{i})=k_{i}^{-1},& S(e_{i})=-k_{i}e_{i},&S(f_{i})=-f_{i}k_{i}^{-1}\\[2ex]
\varepsilon(k_{i})=1, &\varepsilon(e_{i})=0, & \varepsilon(f_{i})=0.
\label{rel8} 
\end{eqnarray}
Letting $q_{h}=e^{h}$ and $K_{i}:=q_{h}^{h_{i}}$ in (\ref{rel1})-(\ref{rel8}), 
we come to definition of $U_{h}(\mathfrak{g})[[h]]$, the topological Hopf algebra over 
$\mathbb{C}[[h]]$.

One introduces a linear ordering on the set of positive roots $\Delta_{+}$ 
by fixing the reduced decomposition of the longest element in the Weyl group
 $w_{0}=s_{i_{1}}s_{i_{2}}\cdots s_{i_{M}}.$ 
Then the linear ordering read from the left to the right is the following  
\begin{equation}
\Delta_{+}=
\{\alpha_{i_{1}}, s_{i_{1}}\alpha_{i_{2}},\ldots, 
s_{i_{1}}s_{i_{2}}\cdots s_{i_{M-1}}\alpha_{i_{M}}\}.
\label{normal}
\end{equation}
An ordering (\ref{normal}) is normal, namely for each $\alpha,\beta\in\Delta_{+}$ 
such that $\alpha+\beta\in\Delta_{+}$ and $\alpha\prec\beta,$ we have 
$\alpha\prec\alpha+\beta\prec\beta.$ There is one-to-one correspondence between the reduced decompositions 
of the longest element in the Weyl group and the normal orderings given by (\ref{normal}).
Following \cite{KT1}, one defines the generators corresponding to the composite roots.  For a chosen
normal ordering on $\Delta_{+}$ let $\alpha,\beta,\gamma\in\Delta_{+}$ be pairwise noncollinear roots,
such that $\gamma=\alpha+\beta.$ Let $\alpha$ and $\beta$ are taken so that there are no other roots 
$\alpha^{\prime}$ and $\beta^{\prime}$ with the property $\gamma=\alpha^{\prime}+\beta^{\prime}.$ Then if 
$e_{\pm\alpha}$ and $e_{\pm\beta}$ have already been constructed, we set
$$ 
\begin{array}{ll}
e_{\gamma}=e_{\alpha}e_{\beta}-q^{-(\alpha,\beta)}e_{\beta}e_{\alpha},&
e_{-\gamma}=e_{-\beta}e_{-\alpha}-q^{-(\beta,\alpha)}e_{-\alpha}e_{-\beta}.
\end{array}
$$ 
  
For any root $\gamma\in\Delta_{+}$ define 
$$
\check{R}_{\gamma}:=
\exp_{q_{\gamma}}(-(q-q^{-1})a_{\gamma}^{-1}\mathop{}e_{\gamma}k_{\gamma}\otimes k_{\gamma}^{-1}f_{\gamma}),
$$
where 
$q_{\gamma}=q^{(\gamma,\gamma)}$  and  
$$
\exp_{q}(x):=\sum_{n\ge 0}\frac{x^{n}}{(n)_{q}!},
\hphantom{aa} (n)_{q}!\equiv (1)_{q}(2)_{q}\ldots (n)_{q},\hphantom{aa} (k)_{q}\equiv (1-q^{k})/(1-q). 
$$
with factors $a_{\gamma}$ coming from the relations
$$
[e_{\gamma}, e_{-\gamma}]=a_{\gamma}\displaystyle\frac{k_{\gamma}-k_{\gamma}^{-1}}{q-q^{-1}}.
$$
The elements $\check{R}_{\gamma}$ are understood to be taken in some 
completion of $U_{q}(\mathfrak{g})\otimes U_{q}(\mathfrak{g})$ 
(the Taylor extension or the $h-$adic \cite{KT1}). 
Now the coproducts of composite root generators can be expressed in terms 
of the adjoint action of the following factors  
$$
\check{R} _{\prec\beta}:=\prod_{\gamma\prec\beta}\check{R} _{\gamma},
$$ 
where the product over all the positive roots such that $\prec\beta$ 
is taken in accordance with the chosen normal ordering. Namely, 
we have the following 
\begin{ccc}[\cite{KT1, KS}]
Consider the canonical isomorphism $\mathfrak{h}\simeq\mathfrak{h}^{*}$ 
defined by the bilinear form 
$(\hphantom{a},\hphantom{a})$
on $\mathfrak{h}.$ Let $h_{\beta}\in\mathfrak{h}$ be the image of a 
root $\beta\in\mathfrak{h}^{*}$
with respect to this isomorphism. Then the following identity holds:
\begin{equation}
\Delta(e_{\beta})=\check{R}_{\prec\beta}(
k_{\beta}^{-1}\otimes e_{\beta}+
e_{\beta}\otimes 1)\check{R}_{\prec\beta}^{-1}.
\label{kt}
\end{equation}
\end{ccc}
\begin{proof}
Note that (\ref{kt}) is equivalent to Proposition 8.3 from \cite{KT1} 
if one applies $(S\otimes S)$ to both parts of 
$$
\Delta^{op}(e_{\alpha})=\left(\prod_{\gamma<\alpha} \check{R}_{\gamma}
\right)
(1\otimes e_{\alpha}+e_{\alpha}\otimes\overline{k}_{\alpha})
\left(\prod_{\gamma<\alpha} \check{R}_{\gamma}\right)^{-1}
$$
(notations are of \cite{KT1}) and check that our construction of the 
modified basis differs by the change $q\leftrightarrow q^{-1}$.
\end{proof}
We conclude this section by adding remarks on specialization \cite{DK}. Let 
${\cal A}=\mathbb{C}[q,q^{-1}]_{(q-1)}$ be a ring of rational 
functions regular at $q=1$ and let   
$\hat{U}_{q}(\mathfrak{g})$ be an ${\cal A}$ subalgebra in $U_{q}(\mathfrak{g})$ 
generated by 
$\left\{e_{i},f_{i},k_{i}^{-1}, \displaystyle\frac{k_{i}-1}{q-1}\right\}$
then 
$$
\hat{U}_{q}(\mathfrak{g})\otimes_{\cal A}\mathbb{C}\approx U(\mathfrak{g})
$$
where $\mathbb{C}$ is regarded as an ${\cal A}$ module ($q$ acts as $1$).
By construction of the quantum Cartan-Weyl basis we have the property 
$$
\Delta(e_{\beta})\in \hat{U}_{q}(\mathfrak{g})\otimes
\hat{U}_{q}(\mathfrak{g})
$$
but, in fact, one can deduce from $(\ref{kt})$ more restrictive property 
\begin{equation}
\Delta(e_{\beta})-k_{\beta}^{-1}\otimes e_{\beta}-
e_{\beta}\otimes 1\in (q-q^{-1})\mathop{}\hat{U}_{q}^{+}(\mathfrak{g})\otimes
\hat{U}_{q}^{+}(\mathfrak{g})
\label{restriction}
\end{equation}
where $\hat{U}_{q}^{+}(\mathfrak{g})$ is generated by 
$\left\{e_{i}, k_{i}^{-1}, \displaystyle\frac{k_{i}-1}{q-1}\right\}$.
In what follows we are usually working rather with completions  
$\hat{U}_{q}(\mathfrak{g})[[t]]$ and 
$\hat{U}_{q}(\mathfrak{g})\otimes\hat{U}_{q}(\mathfrak{g})[[t]]$ in which the 
twists are to be defined. Let us formulate the following simple result which is 
of value for further application 
\begin{aaa}
Let ${\cal F}\in \hat{U}_{q}(\mathfrak{g})
\otimes\hat{U}_{q}(\mathfrak{g})[[t]]$ be a twist in 
$U_{q}(\mathfrak{g})[[t]]$, then its specialization 
$\overline{{\cal F}}$, obtained by order-wise specialization of its 
coefficients from $\hat{U}_{q}(\mathfrak{g})\otimes\hat{U}_{q}(\mathfrak{g})$ at $q=1$, 
is still a twist in $U(\mathfrak{g})[[t]]$.
\label{specialization}
\end{aaa}
\begin{proof}
Indeed,representing ${\cal F}$ as a series 
$$
{\cal F}=1\otimes 1+{\cal F}_{1}t+{\cal F}_{2}t^{2}+\cdots
$$ 
we see that (\ref{drinfeld}) is a equivalent to an infinite set of identities 
and each of them after specializing $q=1$ remain valid. Thus 
$$
\overline{{\cal F}}=1\otimes 1+\overline{{\cal F}}_{1}t+
\overline{{\cal F}}_{2}t^{2}+\cdots
$$    
is a twist in $U(\mathfrak{g})[[t]]$.
\end{proof} 

\section{Quantum Jordanian twists}
We restrict ourselves to consideration of nonexceptional 
Lie algebra $\mathfrak{g}$ and define the quantum Jordanian twists as 
those specializing to semi-classical ones which define quantization of 
skew-symmetric extended Jordanian $r-$matrices: 
$$ 
r_{\lambda}=H_{\lambda}\wedge E_{\lambda}+2\sum_{\gamma_{1}\prec\gamma_{2}, 
\gamma_{1}+\gamma_{2}=\lambda}E_{\gamma_{1}}\wedge 
E_{\gamma_{2}}
$$
by the rule 
$$
{\cal R}={{\cal F}_{\lambda}}_{21}{\cal F}_{\lambda}^{-1}=1\otimes 1+
t\mathop{}r_{\lambda}\bmod{t^{2}}, 
$$
where we have denoted by $H_{\lambda}, E_{\gamma}$ 
the elements of the classical Cartan-Weyl basis.

Let us fix some normal ordering on $\Delta_{+}$ and define 
a generator $e_{\lambda}\in U_{q}(\mathfrak{g})$ corresponding to the 
highest root $\lambda$ according to the recipe from the previous section.
Then nonexceptional root systems are remarkable by the following 
property:
\begin{aaa}
Let $\mathfrak{g}$ be a non exceptional Lie algebra, then there 
is such a normal ordering "$\prec$" on $\Delta_{+}$ so that 
$$
[e_{\gamma},e_{\lambda}]_{q^{-(\gamma,\lambda)}}=0 
\mbox{ for any }\gamma\prec \lambda, \mbox{ and }
e_{\gamma},e_{\lambda}\in U_{q}(\mathfrak{g}).
$$ 
\label{ordering}
\end{aaa}  
\begin{proof}
The proof is based on the expansion \cite{KT1,KS}: 
\begin{equation}
e_{\gamma}e_{\lambda}-q^{-(\gamma,\lambda)}e_{\lambda}e_{\gamma}=
\sum_{\gamma\prec\gamma_{1}\prec\cdots\prec\gamma_{j}\prec\lambda}
c_{l,\gamma,\lambda}e_{\gamma_{1}}^{l_{1}}\cdots e_{\gamma_{j}}^{l_{j}}
\label{ext}
\end{equation}
where $c_{l,\gamma,\lambda}\in\mathbb{C}[q,q^{-1}]$. 
Non zero terms in the sum are subject to condition
\begin{equation}
\gamma+\lambda=l_{1}\gamma_{1}+l_{2}\gamma_{2}+\cdots 
+l_{j}\gamma_{j}.
\label{sss1}
\end{equation}

In $A_{N}$ we choose a normal ordering as 
$$
\alpha_{1}\succ\alpha_{1}
+\alpha_{2}\succ\cdots\succ\alpha_{1}+\alpha_{2}+\cdots+\alpha_{N}\succ\overbrace{\alpha_{2}\succ\cdots\succ\alpha_{N-1}+\alpha_{N}
\succ\alpha_{N}}^{\mbox{roots without $\alpha_{1}$}}
$$ 
and $\lambda=\alpha_{1}+\cdots\alpha_{N}.$
If $\lambda\succ\gamma\succ\alpha_{N}$ we can satisfy (\ref{ext}) only 
with zero coefficients.

In $B_{N}$ we have 
$$
\alpha_{1}\succ\alpha_{1}+\alpha_{2}\succ\cdots\succ\alpha_{1}+2\alpha_{2}+\cdots
+2\alpha_{N-1}+2\beta\succ\overbrace{\alpha_{2}\succ\alpha_{2}+\alpha_{3}\succ\cdots\succ\beta}^{\mbox{roots 
without $\alpha_{1}$}}
$$
and $\lambda=\alpha_{1}+2\alpha_{2}+\cdots
+2\alpha_{N-1}+2\beta$.

In $C_{N}$ we fix the following ordering
$$
\begin{array}{l}
\overbrace{\alpha_{1}\prec\alpha_{1}+\alpha_{2}\prec\cdots\prec\alpha_{1}+\cdots+\alpha_{N-1}}^{\mbox{roots 
without $\beta$}}\prec
2(\alpha_{1}+\cdots+\alpha_{N-1})+\beta\prec\\
 \alpha_{1}+\cdots+\alpha_{N-1}+\beta\prec\cdots\prec
\alpha_{1}+2(\alpha_{2}+\cdots+\alpha_{N-1})+\beta\prec\cdots\prec\alpha_{2}\prec\cdots\prec\beta, 

\end{array}
$$
and $\lambda=2(\alpha_{1}+\cdots+\alpha_{N-1})+\beta$.
This ordering eliminates all non zero terms on the r.h.s of (\ref{ext}).

In $D_{N}$ we have quite a similar situation
$$
\begin{array}{l}
\alpha_{1}\succ\alpha_{1}+\alpha_{2}\succ\cdots\succ\alpha_{1}+\cdots+\alpha_{N-1}\succ 
\alpha_{1}+\cdots+\alpha_{N-1}+\beta
\succ\\
\alpha_{1}+\alpha_{2}+\cdots+2\alpha_{N-2}+\alpha_{N-1}+\beta\succ\cdots\succ
\alpha_{1}+2\alpha_{2}+\cdots 2\alpha_{N-2}
+\alpha_{N-1}+\beta\succ\\[2ex]
\hphantom{aaaaaaaaaaaaaaaaaaaaaaaaaaaaaaaaaaaaaaaaaaaaaaa}\overbrace{\alpha_{2}\succ
\alpha_{2}+\alpha_{3}\succ\cdots\succ\beta}^
{\mbox{roots without $\alpha_{1}$}},
\end{array}
$$ 
and $\lambda=\alpha_{1}+2\alpha_{2}+\cdots +2\alpha_{N-2}
+\alpha_{N-1}+\beta.$
\end{proof}
As a direct consequence of Proposition \ref{ordering}, we obtain 
\begin{aaa}
\label{qcommut}
A $q-$commutation holds 
$$
(e_{\lambda}\otimes 1)(\Delta(e_{\lambda})-e_{\lambda}\otimes 1)=
q_{\lambda}(\Delta(e_{\lambda})-e_{\lambda}\otimes 1)
(e_{\lambda}\otimes 1) \mbox{ where }q_{\lambda}=q^{(\lambda,\lambda)}.
$$ 
\end{aaa} 
\begin{proof}
The statement follows from Proposition \ref{ordering} and (\ref{kt}) 
if one notices that 
$$
\check{R}_{\prec\lambda}(e_{\lambda}\otimes 1)
=(e_{\lambda}\otimes 1)\check{R}_{\prec\lambda}.
$$
\end{proof}
Let us define a coboundary twist
$$
{\cal J}(e_{\lambda})=
(W\otimes W)\Delta(W^{-1}), \mbox{ where }
W=\exp_{q_{\lambda}}(\frac{t}{1-q_{\lambda}}\mathop{}e_{\lambda})
$$  
Then the following is true  
\begin{aaa} 
\label{jordan}
${\cal J}(e_{\lambda})$ is nonsingular and defines a 
nontrivial twisting of $U(\mathfrak{g})[[t]]$ in the limit 
$q\rightarrow 1$. 
\end{aaa}
\begin{proof}
By Proposition \ref{qcommut} we have 
$$
{\cal J}(e_{\lambda})=\exp_{q_{\lambda}}
(\frac{t}{1-q_{\lambda}}\mathop{}1\otimes e_{\lambda})
\exp_{q_{\lambda}^{-1}}(-\frac{t}{1-q_{\lambda}}\mathop{}
(\Delta(e_{\lambda})-e_{\lambda}\otimes 1)).
$$
The latter representation is nonsingular and from 
$\hat{U}_{q}(\mathfrak{g})\otimes\hat{U}_{q}(\mathfrak{g})[[t]$, 
which is obvious if one uses the Campbell-Hausdorff formula after  
applying the dilogarithmic representation of $q-$exponent \cite{FK}:
\begin{equation}
\exp_{q_{\lambda}}(\frac{t}{1-q_{\lambda}}\mathop{}x)=
\exp\left(\sum_{n\ge 1}\frac{t^{n}}{n(1-q_{\lambda}^{n})}\mathop{}x^{n}\right)
\label{dilog}
\end{equation}
along with the properties (\ref{restriction}) and    
$$
[e_{\lambda}, \hat{U}_{q}^{+}(\mathfrak{g})]\in 
(q-1)\mathop{}\hat{U}_{q}^{+}(\mathfrak{g}).
$$
Note, that to be self-consistent one can directly verify that 
(\ref{dilog}) satisfies  (\ref{heine}), considering (\ref{heine}) as 
a functional equation for the function 
${\rm Li}_{2}(t\cdot x,q):$\\ ${\rm Li}_{2}(t\cdot x,q)=
\ln(\exp_{q}(\frac{t}{1-q}\mathop{}x))$. 
Finally, $\overline{{\cal J}}$ is a twist by Proposition \ref{specialization}. 
\end{proof} 
\begin{bbb}
\label{example}
Consider $\mathfrak{g}=\mathfrak{sl}_{N+1}$. We have the following 
formula for the coproduct associated with the chosen normal ordering 
$$
\alpha_{1}\succ\alpha_{1}
+\alpha_{2}\succ\cdots\succ\alpha_{1}+\alpha_{2}+\cdots+\alpha_{N}\succ\overbrace{\alpha_{2}\succ\cdots\succ\alpha_{N-1}+\alpha_{N}
\succ\alpha_{N}}^{\mbox{roots without $\alpha_{1}$}}
$$ 
given by
$$
\Delta(e_{\epsilon_{1}-\epsilon_{N+1}})
=k_{\epsilon_{1}-\epsilon_{N+1}}^{-1}\otimes e_{\epsilon_{1}-\epsilon_{N+1}}+
e_{\epsilon_{1}-\epsilon_{N+1}}\otimes 1
+(1-q^{2})\mathop{}\sum_{i=1}^{N-1}e_{\epsilon_{1}-\epsilon_{i+1}}
k_{\epsilon_{i+1}-\epsilon_{N+1}}^{-1}
\otimes e_{\epsilon_{i+1}-\epsilon_{N+1}}
$$
where $\alpha_{i}=\epsilon_{i}-\epsilon_{i+1}$.

To calculate specialization $\overline{{\cal J}(e_{\epsilon_{1}-\epsilon_{N+1}})}$,
we represent    
$
{\cal J}(e_{\epsilon_{1}-\epsilon_{N+1}})
$
in the following form 
$$
{\cal J}(e_{\epsilon_{1}-\epsilon_{N+1}})=
\exp_{q^{-2}}(-t\sum_{i=1}^{N-1}e_{\epsilon_{1}-\epsilon_{i+1}}
k_{\epsilon_{i+1}-\epsilon_{N+1}}^{-1}
\otimes e_{\epsilon_{i+1}-\epsilon_{N+1}}C_{1,N+1}){\cal J}_{1}
$$ 
where 
\begin{equation}
C_{1,N+1}=\exp_{q^{2}}(-\frac{q\mathop{}t}{1-q^{2}}\mathop{}
e_{\epsilon_{1}-\epsilon_{N+1}})\cdot
\exp_{q^{-2}}(\frac{t}{1-q^{2}}\mathop{}e_{\epsilon_{1}-\epsilon_{N+1}})
\label{quantity1}
\end{equation}
and 
\begin{equation}
{\cal J}_{1}=\exp_{q^{2}}(\frac{t}{1-q^{2}}\mathop{}1\otimes 
e_{\epsilon_{1}-\epsilon_{N+1}})\exp_{q^{-2}}(-\frac{t}{1-q^{2}}\mathop{}
k_{\epsilon_{1}-\epsilon_{N+1}}^{-1}\otimes e_{\epsilon_{1}-\epsilon_{N+1}}).
\label{quantity2}
\end{equation}
Calculation of $(\ref{quantity1})-(\ref{quantity2})$ is based on the 
Heine's formula from \cite{Kac}:
\begin{equation}
1+\sum_{n\ge 1}t^{n}\frac{(\alpha)_{q}\cdots (\alpha+n-1)_{q}}{(n)_{q}!}\mathop{}x^{n}=
\exp_{q}(\frac{t}{1-q}\mathop{}x)
\exp_{q^{-1}}(-\frac{q^{\alpha}t}{1-q}\mathop{}x).
\label{heine}
\end{equation}
Note that (\ref{heine}) can be recast so that to hold in 
$\hat{U}_{q}(\mathfrak{g})\otimes \hat{U}_{q}(\mathfrak{g})[[t]]$: 
\begin{equation}
1\otimes 1+\sum_{n\ge 1}\frac{t^{n}}{(n)_{q^{2}}!}
\left(\frac{k_{\epsilon_{1}-\epsilon_{N+1}}^{-1}-1}{q^{2}-1}\right)\cdots 
\left(\frac{k_{\epsilon_{1}-\epsilon_{N+1}}^{-1}q^{2(n-1)}-1}{q^{2}-1}\right)\mathop{}
\otimes e_{\epsilon_{1}-\epsilon_{N+1}}^{n}={\cal J}_{1}
\label{heine2}
\end{equation}
as one checks 
$\frac{k_{\epsilon_{1}-\epsilon_{N+1}}^{-1}q^{2(n-1)}-1}{q^{2}-1}\in\hat{U}_{q}(\mathfrak{g})$.
Applying the specialization map  
$$
a\mapsto \overline{a\vphantom{f}}:=a\otimes_{\cal A} 1
$$
to each of the tensor factors in 
${\cal J}(e_{\epsilon_{1}-\epsilon_{N+1}})$ we come to a formula of \cite{GZ,KLM}: 
$$
\begin{array}{l}
\overline{{\cal J}(e_{\epsilon_{1}-\epsilon_{N+1}})}=\\[2ex]
\exp(-t\sum_{i=1}^{N-1}E_{1,i+1}\otimes 
E_{i+1,N+1}e^{-\frac 12 \sigma_{1,N+1}})\cdot\\[2ex] 
\displaystyle\left(1\otimes 1+\sum_{n\ge 1}(-1)^{n}t^{n}\frac{H_{1,N+1}(H_{1,N+1}-1)\cdots (H_{1,N+1}-n+1)}
{2^{n}n!}\otimes E_{1,N+1}\right)
\end{array}
$$
where 
$$
\sigma_{1,N+1}=\ln(1-t\mathop{}E_{1,N+1})=
-\sum_{n\ge 1}\frac{t^{n}}{n}\mathop{}E_{1,N+1}
$$
and 
$$
\begin{array}{lcl}
E_{i,j}=\overline{e_{\epsilon_{i}-\epsilon_{j}}\vphantom{l}} && 
H_{1,N}=\overline\frac{k_{\epsilon_{1}-\epsilon_{N+1}}-1\vphantom{l}}{q-1}.
\end{array}
$$
\end{bbb} 
\section{The Cremmer-Gervais twist and its specialization at 
$q\rightarrow 1$}
In this section we consider nontrivial quantum twists in 
$U_{q}(\mathfrak{sl}_{3})[[t]]$ and their semi-classical limits 
$q\rightarrow 1$. As is known 
from the classification of Belavin-Drinfeld triples, there 
are two possible Belavin-Drinfeld triples for $\mathfrak{sl}_{3}$. 
The first one, the empty triple, is accounted for the Drinfeld-Jimbo 
deformation itself, while the second is associated with another 
deformation which can be called the Cremmer-Gervais quantization \cite{CG} 
and there is a solution to (\ref{drinfeld}) defining the twisting element 
providing a possibility to deform 
$U_{h}(\mathfrak{sl}_{3})[[h]]$ further. To be self-consistent we 
first recall the construction of this twist from \cite{KM} and 
then study different possibilities to define specialization 
$q\rightarrow 1$, unveiling a surprising connection with the 
Connes-Moscovici algebra ${\cal H}_{1}$.
\begin{aaa}[\cite{KM}]
\label{cremmer}
An element 
$$
{\cal J}_{CG}=\Phi_{CG}\cdot{\cal K}=\exp_{q_{h}^{-2}}
(\xi\mathop{}e_{32}\otimes e_{12})\mathop{}
q_{h}^{h_{w_{2}}\otimes h_{w_{1}}},\hphantom{aaa}
\xi\in h\cdot\mathbb{C}[[h]],
$$
where 
\begin{equation}
{
\begin{array}{lcl}
h_{w_{1}}=\frac 23 e_{11}-\frac 13 (e_{22}+e_{33}),&&
h_{w_{2}}=\frac 13 (e_{11}+e_{22})-\frac 23 e_{33}
\end{array}
\label{cartans}
} 
\end{equation}
with $w_{1,2}$ being the fundamental weights,
is a twist.
\end{aaa}
\begin{proof}
It is clear that 
${\cal K}=q^{h_{w_{2}}\otimes h_{w_{1}}}_{h}$
defines an abelian twist of $U_{q_{h}}(\mathfrak{sl}_{3})[[h]]$. 
It leads to the following new Hopf algebra  
$U_{q_{h}}^{\cal K}(\mathfrak{sl}_{3})[[h]]$ with the same algebra 
structure as for $U_{q_{h}}(\mathfrak{sl}_{3})[[h]]$ and the new deformed 
coproducts:
$$
\begin{array}{lcl} 
\Delta_{\cal K}(e_{12})=q^{-2\mathop{}h_{1,-1}}_{h}\otimes e_{12}+e_{12}\otimes 1,
&&\Delta_{\cal K}(e_{23})=q^{h_{1,-2}}_{h}\otimes e_{23}+
e_{23}\otimes q^{h_{1,0}}_{h}\\[2ex]
\Delta_{\cal K}(e_{21})=
e_{21}\otimes q_{h}^{h_{2,-1}}+q^{-h_{0,1}}\otimes e_{21},
&&\Delta_{\cal K}(e_{32})=e_{32}
\otimes q^{-2 h_{1,-1}}_{h}+1\otimes e_{32}
\end{array}
$$
where $h_{m,n}:=m\mathop{}h_{w_{1}}+n\mathop{}h_{w_{2}}$.
We are done if we prove that $\Phi_{CG}$ is a twist for 
$U_{q_{h}}(\mathfrak{sl}_{3})[[h]]$. Indeed, explicitly 
(\ref{drinfeld}) reads as the following 
$$
\begin{array}{r}
\exp_{q_{h}^{-2}}(\xi\mathop{}e_{32}\otimes e_{12}\otimes 1)\cdot
\exp_{q_{h}^{-2}}
(\xi\mathop{}(e_{32}
\otimes q^{-2 h_{1,-1}}_{h}+1\otimes e_{32})\otimes e_{12})=\\[2ex]
\exp_{q_{h}^{-2}}(\xi\mathop{}1\otimes e_{32}\otimes e_{12})\cdot
\exp_{q_{h}^{-2}}(\xi\mathop{}e_{32}\otimes
(q^{-2\mathop{}h_{1,-1}}_{h}\otimes e_{12}+e_{12}\otimes 1))
\end{array}
$$
and by the characteristic property of $q-$exponent
$$
\exp_{q}(x+y)=\exp_{q}(y)\exp_{q}(x);\hphantom{aaaa}xy=q\mathop{}yx
$$
(\ref{drinfeld}) holds.  
\end{proof}
Consider now the problem of defining specialization of 
${\cal J}_{CG}$. To do so, we introduce from the beginning a 
$\mathbb{C}(q)$ analogue of $U_{q_{h}}^{\cal K}(\mathfrak{sl}_{3})[[h]]$, 
which we denote by $U^{\prime}_{q}(\mathfrak{sl}_{3})$. 
As an algebra $U^{\prime}_{q}(\mathfrak{sl}_{3})$ is an 
extension of $U_{q}(\mathfrak{sl}_{3})$ obtained by attaching 
elements $L_{i}$ (the maximal lattice), so that 
$K_{j}=\prod_{i=1}^{3}L^{a_{ij}}_{i}$. On the other hand, as a 
coalgebra $U^{\prime}_{q}(\mathfrak{sl}_{3})$ has a new 
coproduct fixed uniquely by its values on the Chevalley generators
$$
\Delta(L_{i})=L_{i}\otimes L_{i}
$$  
$$
\begin{array}{lcl} 
\Delta(e_{1})=
L_{1}^{-2}L_{2}^{2}\otimes e_{1}+e_{1}\otimes 1,
&&\Delta(e_{2})=L_{1}L_{2}^{-2}\otimes e_{2}+
e_{2}\otimes L_{1}\\[2ex]
\Delta(f_{1})=
f_{1}\otimes L_{1}^{2}L_{2}^{-1}+L_{2}^{-1}\otimes f_{1},
&&\Delta(f_{2})=f_{2}
\otimes L_{1}^{-2}L_{2}^{2}+1\otimes f_{2},
\end{array}
$$
where
$L_{i}$ are invertible, pairwise commuting and satisfying  
$$
\begin{array}{lcl}
L_{i}e_{j}L_{i}^{-1}=q^{\delta_{i,j}}e_{j}, &&
L_{i}f_{j}L_{i}^{-1}=q^{-\delta_{i,j}}f_{j}.
\end{array}
$$
In the sense of $\cite{KD}$ define the regular form 
$\hat{U}_{q}^{\prime}(\mathfrak{sl}_{3})$ as a 
${\cal A}:=\mathbb{C}[q,q^{-1}]_{(q-1)}$ algebra such that there is an 
isomorphism  
$$
\hat{U}_{q}^{\prime}(\mathfrak{sl}_{3})
\otimes_{\cal A}\mathbb{C}(q)\approx 
U_{q}^{\prime}(\mathfrak{sl}_{3}).
$$
$\hat{U}_{q}^{\prime}(\mathfrak{sl}_{3})$ is generated over 
${\cal A}$ by
the following set of generators 
$$
\left\{L_{i}^{-1},\displaystyle\frac{L_{i}-1}{q-1}, e_{i}, f_{i}
\right\}
$$ 
If we denote by the "barred" generators the images of generators under 
specialization map 
$$
a\mapsto \overline{a\vphantom{f}}:=a\hat{\otimes}_{\cal A} 1,
$$
then the generators of $\hat{U}_{q}^{\prime}(\mathfrak{sl}_{3})$ specialize 
to the classical Chevalley generators of $U(\mathfrak{sl}_{3})$ as the following:
\begin{equation}
{
\begin{array}{lclcl}
\overline{e}_{1}=E_{12},&& \overline{e}_{2}=E_{23},&&
\displaystyle\overline\frac{L_{1}^{\vphantom{f}}-1}{q-1}=h_{w_{1}}\\[2ex] 
\overline{f}_{1}=E_{21}, && \overline{f}_{2}=E_{32},&&  
\displaystyle\overline\frac{L_{2}^{\vphantom{f}}-1}{q-1}=h_{w_{2}} 
\end{array}
}
\label{limit}
\end{equation}
(see (\ref{cartans})).
Let us additionally make completion $U_{q}^{\prime}(\mathfrak{sl}_{3})[[t]]$ 
then we can formulate a $q-$analogue of Proposition \ref{cremmer}:    
\begin{aaa}
An element 
$$
\hat{\Phi}_{CG}=\exp_{q^{-2}}(\zeta\mathop{}f_{2}
\otimes e_{1}),\hphantom{aaa}
\zeta\in t\cdot{\cal A}[[t]]
$$
is a twist in $\hat{U}_{q}^{\prime}(\mathfrak{sl}_{3})[[t]]$.
\label{crem}
\end{aaa}  
On the other hand apart from $\hat{\Phi}_{CG}$ we can construct  
the twists in $U_{q}^{\prime}(\mathfrak{sl}_{3})[[t]]$ which restrict  
to $\hat{U}_{q}^{\prime}(\mathfrak{sl}_{3})[[t]]$ only after a suitable  
change of basis which is implemented by some coboundary twist, 
namely we can prove 
\begin{aaa}
An element 
$$
{\cal F}_{CG}=
(V\otimes V)\exp_{q^{-2}}(\frac{q^{-2}\cdot\zeta\cdot t}{1-q^{2}}\mathop{}f_{2}\otimes e_{1})
\Delta(V^{-1}),
\hphantom{aaa}
\zeta\in t\cdot{\cal A}[[t]],
$$
where
$$
V=
\exp_{q^{-2}}(-\frac{\zeta}{1-q^{2}}\mathop{}
e_{1}L_{1}^{2}L_{2}^{-2})\cdot\exp_{q^{2}}(\frac{t}{1-q^{2}}\mathop{}f_{2}),
$$
restricts to a twist of $\hat{U}_{q}^{\prime}(\mathfrak{sl}_{3})[[t]]$.
\begin{proof}
By the form of the coproducts 
$$
\begin{array}{lcl}
\Delta(e_{1}L_{1}^{2}L_{2}^{-2})=e_{1}L_{1}^{2}L_{2}^{-2}\otimes 
L_{1}^{2}L_{2}^{-2}+1\otimes e_{1}L_{1}^{2}L_{2}^{-2},&& 
\Delta(f_{2})=
f_{2}\otimes L_{1}^{-2}L_{2}^{2}+1\otimes f_{2}
\end{array}
$$
we have explicitly 
$$
\begin{array}{l}
{\cal F}_{CG}=\left(
V\otimes \exp_{q^{-2}}(-\displaystyle\frac{\zeta}{1-q^{2}}\mathop{}e_{1}L_{1}^{2}L_{2}^{-2})
\right)\\[2ex]
\displaystyle\exp_{q^{-2}}(\frac{q^{-2}\cdot\zeta\cdot t}{1-q^{2}}\mathop{}
f_{2}\otimes e_{1})\exp_{q^{-2}}(-\frac{t}{1-q^{2}}\mathop{}
f_{2}\otimes L_{1}^{-2}L_{2}^{2})\Delta(
\exp_{q^{2}}(\frac{\zeta}{1-q^{2}}\mathop{}e_{1}L_{1}^{2}L_{2}^{-2})).
\end{array}
$$
Using the five terms relation, \cite{FK}:\\
If $[u,[u,v]]_{q^{2}}=[v,[u,v]]_{q^{-2}}=0$ then
\begin{equation}
{\rm e}_{q^{2}}(u)\cdot{\rm e}_{q^{2}}(v)={\rm e}_{q^{2}}(v)\cdot
{\rm e}_{q^{2}}(\displaystyle\frac{1}{1-q^{2}}\mathop{}[u,v])\cdot{\rm e}_{q^{2}}(u)
\end{equation}
we can simplify
$$
\begin{array}{l}
\displaystyle\exp_{q^{-2}}(-\frac{t}{1-q^{2}}\mathop{}f_{2}\otimes 
L_{1}^{-2}L_{2}^{2})\exp_{q^{-2}}(-\frac{\zeta}{1-q^{2}}\mathop{}
1\otimes e_{1}L_{1}^{2}L_{2}^{-2})=\\[2ex]
\displaystyle\exp_{q^{-2}}(-\frac{\zeta}{1-q^{2}}\mathop{}
1\otimes e_{1}L_{1}^{2}L_{2}^{-2})
\exp_{q^{-2}}(\frac{q^{-2}\cdot\zeta\cdot t}{1-q^{2}}\mathop{}
f_{2}\otimes e_{1})\exp_{q^{-2}}(-\frac{t}{1-q^{2}}\mathop{}f_{2}\otimes 
L_{1}^{-2}L_{2}^{2})
\end{array}
$$
and thus 
${\cal F}_{CG}$ transforms to the following form  
$$
{\cal F}_{CG}=\displaystyle\exp_{q^{-2}}(-\frac{t}{1-q^{2}}\mathop{}(f_{2}\otimes 
L_{1}^{-2}L_{2}^{2}+e_{1}L_{1}^{2}L_{2}^{-2}\otimes 1))
\exp_{q^{2}}(\frac{t}{1-q^{2}}\mathop{}(f_{2}\otimes 1+
L_{1}^{2}L_{2}^{-2}e_{1}\otimes L_{1}^{2}L_{2}^{-2}))
$$
and finally by the Heine's formula we have  
$$
{\cal F}_{CG}=1\otimes 1+\sum_{n\ge 1}\frac{1}{(n)_{q^{2}}!}
\left(1\otimes \frac{L_{1}^{-2}L_{2}^{2}-1}{q^{2}-1}\cdots 
\frac{L_{1}^{-2}L_{2}^{2}q^{2(n-1)}-1}{q^{2}-1}\right)(t\cdot f_{2}\otimes 1+
\zeta\cdot L_{1}^{2}L_{2}^{-2}e_{1}\otimes L_{1}^{2}L_{2}^{-2})^{n}
$$
the latter expression restricts to 
$\hat{U}_{q}^{\prime}(\mathfrak{sl}_{3})[[t]]$ 
as we  have 
$$
\frac{L_{1}^{-2}L_{2}^{2}q^{2(n-1)}-1}{q^{2}-1}\in 
\hat{U}_{q}^{\prime}(\mathfrak{sl}_{3})
$$
\end{proof}
\end{aaa}
\section{Semi-classical twists and noncommutative geometry}
\begin{aaa}
There is a semi-classical twist ${\cal F}$ and a homomorphism 
\label{embedding}
$$
\iota:
{\cal H}_{1}\rightarrow U^{\cal F}(\mathfrak{sl}_{3})[[t]]
$$
\end{aaa}
\begin{proof}
We solve more general problem of obtaining quantization 
${\cal H}_{1,q}^{\prime}$ in the sense that there is an embedding 
$$  
\iota_{q}: {\cal H}_{1,q}^{\prime}\hookrightarrow 
U_{q}^{\prime}(\mathfrak{sl}_{3})[[t]]
$$
where ${\cal H}_{1,q}^{\prime}$ denotes 
an appropriately defined $q-$deformation of 
${\cal H}_{1}^{\prime}$. 
Consider the following coboundary 
twist in $U_{q}^{\prime}(\mathfrak{sl}_{3})[[t]]$:   
$$
{\cal J}:=(W\otimes W)\Delta(W^{-1}),
\hphantom{aaa}W=\exp_{q^{-2}}(-\frac{t}{1-q^{2}}\mathop{}e_{1+2}L_{1})
$$
where
$$
e_{1+2}L_{1}:=(e_{1}e_{2}-q\mathop{}e_{2}e_{1})L_{1}
$$
has the following coproduct
$$
\Delta(e_{1+2}L_{1})=e_{1+2}L_{1}\otimes L_{1}^{2}+
1\otimes e_{1+2}L_{1}+(1-q^{2})\mathop{}e_{1}L_{1}^{2}L_{2}^{-2}
\otimes e_{2}L_{1}. 
$$
By the properties 
$$
\begin{array}{l}
(1\otimes e_{1+2}L_{1})(e_{1+2}L_{1}\otimes L_{1}^{2}+
(1-q^{2})\mathop{}e_{1}L_{1}^{2}L_{2}^{-2}\otimes e_{2}L_{1})=
\\[2ex]
\hphantom{aaaaaaaaaaaaaa}q^{-2}
(e_{1+2}L_{1}\otimes L_{1}^{2}+
(1-q^{2})e_{1}L_{1}^{2}L_{2}^{-2}\otimes e_{2}L_{1})
(1\otimes e_{1+2}L_{1}),\\[2ex]
(e_{1}L_{1}^{2}L_{2}^{-2}\otimes e_{2}L_{1})
(e_{1+2}L_{1}\otimes L_{1}^{2})
=q^{-2}(e_{1+2}L_{1}\otimes L_{1}^{2})
(e_{1}L_{1}^{2}L_{2}^{-2}\otimes e_{2}L_{1})
\end{array}
$$
${\cal J}$ is nonsingular, the reasoning is 
same as in Proposition \ref{jordan}, and can be represented 
in the following form
$$
\begin{array}{r}
{\cal J}=
\overbrace{\displaystyle{\rm Ad}\exp_{q^{-2}}(-\frac{t}{1-q^{2}}\mathop{}e_{1+2}L_{1}\otimes 1)
(\exp_{q^{2}}(t\mathop{}e_{1}L_{1}^{2}L_{2}^{-2}\otimes e_{2}L_{1}))}^{{\cal J}_{1}}\cdot
\\[2ex]
\overbrace
{\displaystyle\exp_{q^{-2}}(-\frac{t}{1-q^{2}}\mathop{}e_{1+2}L_{1}\otimes 1)
\exp_{q^{2}}(\frac{t}{1-q^{2}}\mathop{}
e_{1+2}L_{1}\otimes L_{1}^{2})}^{{\cal J}_{2}}.
\end{array}
$$
Let us check that for both factors we have 
$$
{\cal J}_{1,2}\in
 \hat{U}_{q}^{\prime}(\mathfrak{sl}_{3})\otimes 
 \hat{U}_{q}^{\prime}(\mathfrak{sl}_{3})[[t]].
 $$
Indeed, it follows from explicit form of the factors 
${\cal J}_{1,2}$:
$$
{\cal J}_{1}=
\exp_{q^{2}}\left(t\mathop{}
e_{1}L_{1}^{2}L_{2}^{-2}\mathop{}\frac{1}{1-t\mathop{}e_{1+2}L_{1}}\otimes 
e_{2}L_{1}\right),
$$
$$
\begin{array}{l}
{\cal J}_{2}=
\displaystyle\exp_{q^{-2}}\left(\frac{q^{-2}\mathop{}t}{1-q^{-2}}\mathop{}
e_{1+2}L_{1}\otimes 1\right)
\left(\exp_{q^{-2}}\left(\frac{q^{-2}\mathop{}t}{1-q^{-2}}\mathop{}
e_{1+2}L_{1}\otimes L_{1}^{2}\right)\right)^{-1}=\\[2ex]
=\displaystyle 1\otimes 1+\sum_{n\ge 1}t^{n}\frac{(-1)^{n}}{(n)_{q^{-2}}!}
(e_{1+2}L_{1})^{n}\otimes\frac{(L_{1}^{2}-1)}{q^{2}-1}\cdot
\frac{(L_{1}^{2}q^{-2}-1)}{q^{2}-1}\cdots
\frac{(L_{1}^{2}q^{-2(n-1)}-1)}{q^{2}-1}.
\end{array}
$$
Calculating specialization at $q=1$ we come to 
$$ 
\overline{{\cal J}}=
\exp(t\mathop{}E_{12}\frac{1}{1-t\mathop{}E_{13}}
\otimes E_{23})\left(1\otimes 1+\sum_{n\ge 1}
(-1)^{n}t^{n}E_{13}^{n}\otimes 
\frac{h_{w_{1}}(h_{w_{1}}-1)\cdots (h_{w_{1}}-n+1)}{n!}\right) 
$$
$\overline{\cal J}$ defines a noncoboundary deformation of 
$U(\mathfrak{sl}_{3})[[t]]$ as it follows from 
$$
\overline{\cal J}_{21}\ne \overline{\cal J}.
$$
On the other hand its quantum counter part ${\cal J}$ is coboundary in  
$U_{q}^{\prime}(\mathfrak{sl}_{3})[[t]]$ and amounts to 
switching to another basis of $U_{q}^{\prime}(\mathfrak{sl}_{3})[[t]]$.
In particular, the subset of generators 
$$
\{L_{1}L_{2}^{-1}, e_{1},f_{2}\}
$$
is changing in the following way 
$$
\begin{array}{lcl}
\displaystyle L_{1}^{2}L_{2}^{-2}\mapsto 
W (L_{1}^{2}L_{2}^{-2})W^{-1}=L_{1}^{2}L_{2}^{-2},&&
\displaystyle e_{1}\mapsto W 
(e_{1})W^{-1}=e_{1}\frac{1}{1-t\mathop{}e_{1+2}L_{1}}\\[2ex]
\end{array}
$$
$$
\displaystyle 
f_{2}\mapsto W(f_{2})W^{-1}=f_{2}-
\frac{t}{1-q^{2}}\mathop{}q^{-1}
L_{1}^{2}L_{2}^{-2}e_{1}\mathop{}\frac{1}{1-t\mathop{}
e_{1+2}L_{1}}
$$
and we can form a Hopf subalgebra 
${\cal D}_{q}\subset U_{q}^{\prime {\cal J}}(\mathfrak{sl}_{3})[[t]]$ 
generated by the set of generators  
$$ 
\left\{k:=L_{1}^{2}L_{2}^{-2},x:=f_{2},
z:=q^{-1}\mathop{}L_{1}^{2}L_{2}^{-2}e_{1}\frac{1}{1-t\mathop{}e_{1+2}L_{1}}
\right\}.
$$
with the following structure:
$$
\begin{array}{lclcl}
k xk^{-1}=q^{2}\mathop{}x,&& k z k^{-1}=q^{2}\mathop{}z,
&&
q^{2}xz-zx=-t z^{2}\\[2ex]
\end{array}
$$ 
$$
\begin{array}{lcl}
\Delta(k)=k\otimes k, &&
\Delta(z)=z\otimes k+1\otimes z
\displaystyle
\end{array}
$$
$$
\Delta(x)=x\otimes k^{-1}+1\otimes x+t\mathop{}z\otimes 
\frac{(k-k^{-1})}{1-q^{2}}.
$$
The structure of its specialization ${\cal D}_{1}\approx 
{\cal D}_{q}\hat{\otimes}_{\cal A}1$ is the following:  
$$
\begin{array}{lclcl}
[\overline{\vphantom{f}y},\overline{\vphantom{f}x}]=\overline{\vphantom{f}x},&& 
[\overline{\vphantom{f}y},\overline{\vphantom{f}z}]=\overline{\vphantom{f}z},
&&
\overline{\vphantom{f}x}\overline{\vphantom{f}z}-
\overline{\vphantom{f}z}\overline{\vphantom{f}x}=
-t \overline{\vphantom{f}z}^{2}\\[2ex]
\end{array}
$$ 
$$
\begin{array}{lcl}
\Delta(\overline{\vphantom{f}y})=
\overline{\vphantom{f}y}\otimes 1+1\otimes\overline{\vphantom{f}y}, &&
\Delta(\overline{\vphantom{f}z})=\overline{\vphantom{f}z}\otimes 1+
1\otimes\overline{\vphantom{f}z}
\end{array}
$$
$$
\Delta(\overline{\vphantom{f}x})=\overline{\vphantom{f}x}
\otimes 1+1\otimes\overline{\vphantom{f}x}-t\mathop{}\overline{\vphantom{f}z}\otimes\overline{\vphantom{f}y}
$$
where 
$$
\overline{\vphantom{f}y}:=\overline\frac{\vphantom{l}k-1}{q-1}.
$$
Finally we obtain the stated map
$$  
\iota: {\cal H}_{1}^{\prime}\rightarrow{\cal  D}_{1}\subset 
U^{\overline{\cal J}}(\mathfrak{sl}_{3})[[t]]
$$
by fixing its values on the generators 
$$
\begin{array}{lclcl}
\iota(X)=-\frac 12\overline{\vphantom{f}x}, &&  
\iota(Y)=\overline{\vphantom{f}y}, &&
\iota(Z)=t\overline{\vphantom{f}z}.
\end{array}
$$
where the generators $X,Y,Z$ fulfills the relations 
$$
\begin{array}{lclcl}
[Y,X]=X, && [Y,Z]=Z, && [X,Z]=\frac 12 Z^{2}
\end{array}
$$
and the coproducts are obtained from (\ref{coproducts}) by substitution 
$Z$ for $\delta_{1}$.
Next, $\iota$ is an isomorphism as $\iota$ maps the basis of 
${\cal H}_{1}^{\prime}:$  
$\left\{X^{k}Y^{l}Z^{m}\right\}_{k,l,m\ge 0}$  
onto the basis of ${\cal D}_{1}:$ $\left\{\overline{\vphantom{f}x}^{k}
\overline{\vphantom{f}y}^{l}\overline{\vphantom{f}z}^{m}\right\}_{k,l,m\ge 0}$ 
\end{proof}
Now we can obtain  
\begin{aaa} 
An element 
$$
F_{q}=
(WV\otimes WV)\exp_{q^{-2}}(\frac{q^{-3}\cdot t^{3}}{(1-q^{2})^{2}}\mathop{}f_{2}\otimes e_{1})
\Delta(V^{-1})(W^{-1}\otimes W^{-1}),
$$
where
$$
V=
\exp_{q^{-2}}(-\frac{q^{-1}t}{(1-q^{2})^{2}}\mathop{}
e_{1}L_{1}^{2}L_{2}^{-2})\cdot\exp_{q^{2}}(\frac{t}{1-q^{2}}\mathop{}f_{2}),
$$
restricts to a twist in $\hat{U}_{q}^{\prime {\cal J}}(\mathfrak{sl}_{3})[[t]]$.
\end{aaa}
\begin{proof}
It is convenient to introduce notation  
$$
\tilde{{\cal F}}_{CG}=(V\otimes V)
\exp_{q^{-2}}(\frac{q^{-3}\cdot t^{3}}{(1-q^{2})^{2}}\mathop{}
f_{2}\otimes e_{1})\Delta(V^{-1}).
$$
Then by the reasoning of Proposition \ref{crem}
we can prove that 
$$
\tilde{{\cal F}}_{CG}=1\otimes 1+\sum_{n\ge 1}\frac{1}{(n)_{q^{2}}!}
\left(1\otimes \frac{L_{1}^{-2}L_{2}^{2}-1}{q^{2}-1}\cdots 
\frac{L_{1}^{-2}L_{2}^{2}q^{2(n-1)}-1}{q^{2}-1}\right)(t\cdot f_{2}\otimes 1+
\frac{q^{-1}t^{2}}{1-q^{2}}\cdot L_{1}^{2}L_{2}^{-2}e_{1}\otimes L_{1}^{2}L_{2}^{-2})^{n}
$$
Next, the following element 
$$
\tilde{{\cal F}}_{CG}^{W}:=(W\otimes W)\tilde{{\cal F}}_{CG}\Delta(W^{-1})
$$
is a twist in $U_{q}^{\prime}(\mathfrak{sl}_{3})[[t]]$ and respectively 
$$
F_{q}=(W\otimes W)\tilde{{\cal F}}_{CG}(W^{-1}\otimes W^{-1})
$$
defines a twist in $U_{q}^{\prime{\cal J}}(\mathfrak{sl}_{3})[[t]]$. 
Explicitly
$$ 
\begin{array}{r}
F_{q}=\displaystyle 1\otimes 1+\sum_{n\ge 1}\frac{1}{(n)_{q^{2}}!}
\left(1\otimes \frac{L_{1}^{-2}L_{2}^{2}-1}{q^{2}-1}\cdots 
\frac{L_{1}^{-2}L_{2}^{2}q^{2(n-1)}-1}{q^{2}-1}\right)\cdot\\[2ex]
\displaystyle\left(t\cdot f_{2}\otimes 1-
t^{2}\cdot q^{-1}\mathop{}e_{1}\frac{1}{1-t\mathop{}e_{1+2}L_{1}}\otimes 
\frac{L_{1}^{2}L_{2}^{-2}-1}{q^{2}-1}\right)^{n}
\end{array}
$$
again similarly to Proposition \ref{crem} we check that 
$F_{q}$ resticts to a twist in $\hat{U}_{q}^{\prime {\cal J}}(\mathfrak{sl}_{3})[[t]]$.
If we specialize $q=1$ and apply Proposition \ref{embedding} we obtain 
$$
F_{1}=1\otimes 1+\sum_{n\ge 1}\frac{t^{n}}{n!}
(1\otimes Y(Y-1)\cdots(Y-n+1))(2X\otimes 1+Z\otimes Y)^{n}.
$$
Note that a formula for $F_{1}$ was obtained in \cite{KST} by a direct 
check in the study of the twists for $\mathfrak{sl}_{2}$ Yangian.  
To make correspondence with (\ref{connestwist}) we additionally twist 
$F_{1}$ by a coboundary:
$$
F_{1}^{\prime}=\left(\exp(t\mathop{}XY)\otimes \exp(t\mathop{}XY)\right)F_{1}
\Delta(\exp(-t\mathop{}XY))
$$
expanding in $t$ we have   
$$
F_{1}^{\prime}=1\otimes 1+t\cdot (X\otimes Y-Y\otimes X+ZY\otimes Y)+\cdots
$$ 
Thus $F_{1}^{\prime}$ is equivalent to (\ref{connestwist}).
\end{proof}

\end{document}